\documentclass[11pt]{amsart}
\renewcommand{\bar}{\overline}
\usepackage{amsmath,amsthm,amscd,euscript}
\renewcommand{\bar}{\overline}
\def \r{\mathbb R}

\def \z{\mathbb Z}

\DeclareMathOperator{\MOD}{mod}

\newtheorem{theorem}{Theorem}[section]
\newtheorem{lemma}[theorem]{Lemma}
\newtheorem{statement}[theorem]{Statement}
\newtheorem{proposition}[theorem]{Proposition}
\newtheorem{corollary}[theorem]{Corollary}
\theoremstyle{remark}
\newtheorem{remark}[theorem]{Remark}
\theoremstyle{definition}
\newtheorem{definition}[theorem]{Definition}

\newtheorem{problem}{Problem}

\title{Classification of lattice-regular lattice convex polytopes.}
\author{Oleg~Karpenkov}
\date{13 March 2006}
\thanks{{\it MSC-class}: 11H06; 51M20}
\thanks{
Partially supported by RFBR grant SS-1972.2003.1, by RFBR grant
05-01-01012a, by RFBR 05-01-02805-CNRSL\_a grant, and by NWO-RFBR
grant 047.011.2004.026 (RFBR 05-02-89000-NWO\_a)}


\email[Oleg Karpenkov]{karpenk@mccme.ru}

\address{CEREMADE - UMR 7534
Universit\'e Paris-Dauphine France -- 75775 Paris SEDEX 16}

\begin{document}
\input epsf

\begin{abstract}
In this paper for any dimension $n$ we give a complete description
of lattice convex polytopes in $\r^n$ that are regular with
respect to the group of affine transformations preserving the
lattice.
\end{abstract}

 \maketitle

\tableofcontents

\sloppy \normalsize

\section*{Introduction.}

Consider an $n$-dimensional real vector space. Let us fix a
full-rank lattice in it. A {\it convex polytope} is a convex hull
of a finite number of points. A hyperplane $\pi$ is said to be
{\it supporting} for a (closed) convex polytope $P$, if the
intersections of $P$ and $\pi$ is not empty, and the whole
polytope $P$ is contained in one of the closed half-spaces
bounded by $\pi$. An intersection of any polytope $P$ with any
its supporting hyperplane is called a {\it face} of the polytope.
Zero- and one-dimensional faces are called {\it vertices} and
{\it faces}.

Consider an arbitrary $n$-dimensional convex polytope $P$. An
arbitrary unordered $(n{+}1)$-tuple of faces containing the whole
polytope $P$, some its hyperface, some hyperface of this
hyperface, and so on (up to a vertex of $P$) is called a {\it
face-flag} for the polytope $P$.

A convex polytope is said to be lattice if all its vertices are
lattice points. An affine transformation is called {\it
lattice-affine} if it preserves the lattice. Two convex lattice
polytopes are said to be {\it lattice-congruent} if there exist a
lattice-affine transformation taking one polytope to the other. A
lattice polytope is called {\it lattice-regular} if for any two
its face-flags there exist a lattice-affine transformation
preserving the polytope and taking one face-flag to the other.

In this paper we give a complete description of lattice-regular
convex lattice polytopes in $\r^n$ for an arbitrary $n$
(Theorem~\ref{main_theorem} in Section~\ref{formulation}).

\vspace{2mm}

The study of convex lattice polytopes is actual in different
branches of mathematics, such as lattice geometry(see, for
example~\cite{Bar1}, \cite{Bar2}, \cite{JMK1}, \cite{Rez}),
geometry of toric varieties (see~\cite{Ful}, \cite{Kh1},
\cite{Oda}) and multidimensional continued fractions
(see~\cite{Arn2}, \cite{Kor2}, \cite{KarPyr}, \cite{LacBook},
\cite{Mou2}). Mostly, it is naturally to study such polytopes
with respect to the lattice-congruence equivalence relation.

Now we formulate two classical examples of unsolved problems on
convex lattice polytopes. The first one comes from the geometry
of toric varieties.

\begin{problem}
Find a complete invariant of lattice-congruence classes of convex
lattice (two-dimensional) polygons.
\end{problem}

Only some estimates are known at this moment (see for
example~\cite{Arn5} and~\cite{Bar2}).

The second problem comes from lattice geometry and theory of
multi-dimensional continued fractions. A lattice symplex is
called {\it empty} if the intersection of this (solid) symplex
with the lattice coincides with the set of its vertices.

\begin{problem}
Find a description of lattice-congruence classes of empty
symplices.
\end{problem}

The answer to the second problem in the two-dimensional case is
simple. All empty triangles are lattice-congruent.
Tree-dimensional case is much more complicated. The key to the
description gives White's theorem (1964) shown in~\cite{Whi} (for
more information see~\cite{Rez}, \cite{Mou2}, and~\cite{KarPyr}).

The problems similar to the shown above are complicated and seem
not to be solved in the nearest future. Nevertheless, specialists
of algebraic geometry or theory of multidimensional continued
fractions usually do not need the complete classifications but
just some special examples.

In the present paper we make the first steps in the study of the
lattice polytopes with non-trivial group of {\it
lattice-symmetries} (i.e. the group of lattice-affine
transformations, preserving the polytope). We describe the
``maximally'' possible lattice-symmetric polytopes: the
lattice-regular polytopes.

Let us formulate statement for the second step in the study of
the lattice polytopes with non-trivial group of
lattice-symmetries. A convex lattice pyramid $P$ with the base
$B$ is said to be {\it lattice-regular} if $B$ is a
lattice-regular polytope, and the group of lattice-symmetries of
the base $B$ (in the hyperplane containing $B$) is expandable to
the group of lattice-symmetries of the whole pyramid $P$.

\begin{problem}
Find a description of lattice-regular convex lattice pyramids.
\end{problem}

This paper is organized as follows. We give a well-known
classical description of Euclidean and abstract regular polytopes
in terms of Schl\"afli symbols in Section 1. In Section 2 we give
necessary definitions of lattice geometry and formulate a new
theorem on lattice-affine classification of lattice-regular
convex polytopes. Further in Section~2 we prove this theorem for
the two-dimensional case. We study the cases (in any dimension)
of lattice-regular symplices, cubes, and generalized octahedra in
Sections~4, 5, and 6 respectively. Finally in Section~7 we
investigate the remaining cases of low-dimensional polytopes and
conclude the proof of the main theorem.

{\bf Acknowledgement.} The author is grateful to professors
V.~I.~Arnold and A.~G.~Khovanskii for useful remarks and
discussions, and Universit\'e Paris-Dauphine --- CEREMADE  for
the hospitality and excellent working conditions.

\section{Euclidean and abstract regular polytopes.}

For the proof of the main theorem on lattice-regular polytopes we
use the classification of abstract convex polytopes. We start
this section with the description of Euclidean regular polytopes
and Schl\"afli symbols for them, and then continue with the case
of abstract regular polytopes.

\subsection{Euclidean regular polytopes.}

\begin{figure}[h]
$$\epsfbox{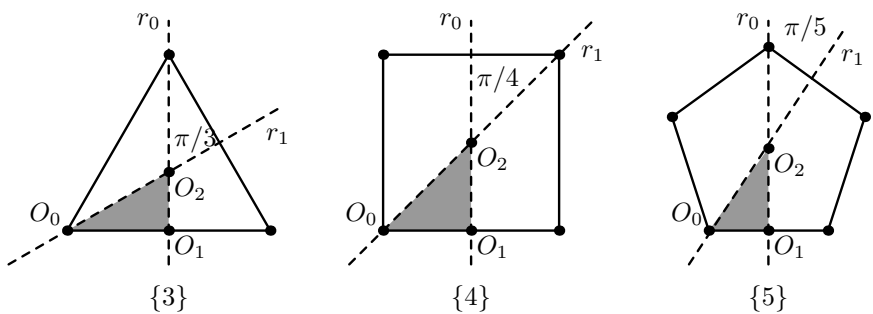}$$
\caption{Basic reflections and Schl\"afli symbols for some
regular polygons.}\label{SchlSymb}
\end{figure}

Consider an arbitrary $n$-dimensional Euclidean regular polytope
$P$. Let $(F_n{=}P,F_{n-1}, \ldots, F_1,F_0)$ be one of its flags.
Denote by $O_i$ the mass center of the face $F_i$ considered as a
homogeneous solid body (for $i=0, \ldots, n$). The
$n$-dimensional tetrahedron $O_{0}O_{1}\ldots O_{n-1}O_{n}$ is
called the {\it chamber} of a regular polytope $P$ corresponding
to the given flag. Denote by $r_i$ (for $i=0, \ldots, {n-1}$) the
reflection about the $(n{-}1)$-dimensional plane spanning the
points $F_n,\ldots, F_{i+1}$, $F_{i-1},\ldots,F_0$. See
Figure~\ref{SchlSymb}. These reflections are sometimes called
{\it basic}.

The following classical statement holds.
\begin{statement}
The reflections $r_0, r_1,\ldots, r_{n-1}$ generate
the group of Euclidean symmetries of the Euclidean regular polytope $P$.\\
For $i=1,\ldots, n{-}1$ the angle between the fixed hyperplanes
of the symmetries $r_{i-1}$ and $r_i$ equals $\pi /a_i$, where
$a_i$ is an integer greater than or equivalent to $3$.
\end{statement}

The symbol $\{a_1,\ldots, a_{n-1} \}$ is said to be the {\it
Schl\"afli symbol} for the polygon $P$. Traditionally, the string
$a,a, \ldots, a$ of the length $s$ in Schl\"afli symbol is
replaced by the symbol $a^s$.

Since all face-flags of any regular polytope are congruent, the
Schl\"afli symbol is well-defined.

\vspace{2mm}

{\bf Theorem A. On classification of regular Euclidean polytopes.}
{\it Any regular convex Euclidean polytope is homothetic to some
polytope of the following list.

{\bf List of regular Euclidean polytopes.}\\
{\bf Dimension 1:} a segment with Schl\"afli symbol $\{\}$.\\
{\bf Dimension 2:} a regular polygon with $m$ vertices $($for any
$m\ge 3$$)$
with Schl\"afli symbol $\{m\}$.\\
{\bf Dimension 3:} a regular tetrahedron $(\{3,3\})$, a regular
octahedron $(\{3,4\})$, a regular cube $(\{4,3\})$, a regular
icosahedron $(\{3,5\})$,
a regular dodecahedron $(\{5,3\})$.\\
{\bf Dimension 4:} a regular symplex $(\{3,3,3\})$, a regular
cube $(\{4,3,3\})$, a regular generalized octahedron $($or cross
polytope, or hyperoctahedron; with Schl\"afli symbol
$\{3,3,4\}$$)$, a regular 24-cell $($or hyperdiamond, or
icositetrachoron; with Schl\"afli symbol $\{3,4,3\}$$)$, a
regular 600-cell $($or hypericosahedron, or hexacosichoron; with
Schl\"afli symbol $\{3,3,5\}$$)$, a regular 120-cell $($or
hyperdodecahedron, or hecatonicosachoronor; with Schl\"afli
symbol $\{5,3,3\}$$)$. {\bf Dimension n (n$>$4):} a regular
symplex $(\{3^{n-1}\})$, a regular cube $(\{4,3^{n-2}\})$, a
regular generalized octahedra $(\{3^{n-2},4\})$. } \vspace{2mm}

\begin{remark}
The cases of dimension one, two, and three were already known to
the ancient mathematicians. The cases of higher dimensions were
studied by Schl\"afli (see in~\cite{Sch}).
\end{remark}

\subsection{Abstract regular polytopes.}

In this subsection we consider arbitrary convex polytopes.
Consider two $n$-dimensional polytopes. A homeomorphism of two
$n$-dimensional polytopes is said to be {\it combinatorical} if
it takes any face of one polytope to some face of the same
dimension of the other polytope. Two polytopes are called {\it
combinatorically isomorphic} if there exist a combinatorical
homeomorphism between them.

A convex polytope is called {\it combinatorical regular} if for
any two its face-flags there exist a combinatorical homeomorphism
taking the polytope to itself and one face-flag to the other.

{\bf Theorem B. (McMullen~\cite{McM1}.)} {\it A polytope is
combinatorical regular iff it is combinatorically isomorphic to a
regular polytope.}

The proof of this statement essentially uses the work of
Coxeter~\cite{Cox}.

\begin{remark}
Theorem~B implies the classification of real affine and
projective polytopes (see~\cite{McM2}). Both classifications
coincide with the classification of Euclidean regular polytopes.
For further investigations of abstract polytopes see for example
the work of L.~Danzer and E.~Schulte~\cite{Dan} and the book on
abstract regular polytopes by P.~McMullen and
E.~Schulte~\cite{McM3}.
\end{remark}

\section{Definitions and formulation of the main result.}\label{formulation}

Let us fix some basis of lattice vectors $\bar e_i$ for
$i=1,\ldots,n$ generating the lattice in $\r^n$. Denote by $O$
the origin in $\r^n$.

Consider arbitrary non-zero integers $n_1, \ldots, n_k$ for $k\ge
2$. By $\gcd(n_1,\ldots,n_k)$ we denote the greater common divisor
of the integers $n_i$, where $i=1,\ldots, k$. We write that
$a\equiv b (\MOD c)$ if the reminders of $a$ and $b$ modulo $c$
coincide.

\subsection{Some definitions of lattice geometry.}
Let $Q$ be an arbitrary lattice polytope with the vertices
$A_i=O+\bar v_i$ (where $\bar v_i$ --- lattice vectors) for
$i=1,\ldots, m$, and $t$ be an arbitrary positive integer. The
polygon $P$ with the vertices $B_i=O+t\bar v_i$ for $i=1,\ldots,
m$ is said to be the {\it $t$-multiple} of the polygon $Q$.

\begin{definition}
A lattice polytope $P$ is said to be {\it elementary} if for any
integer $t>1$ and any lattice polytope $Q$ the polytope $P$ is
not lattice-congruent to the $t$-multiple of the lattice polytope
$Q$.
\end{definition}

\subsection{Notation for particular lattice polytopes.}

We will use the following notation.

{\bf Symplices.} For any $n>1$ we denote by $\{3^{n-1}\}^{L}_p$
the $n$-dimensional symplex with the vertices:
$$
V_0=O, \qquad \hbox{$V_i=O{+}\bar e_i$, for $i=1,\ldots,n{-}1$},
\quad \hbox{and} \quad V_n=(p{-}1)\sum\limits_{k{=}1}^{n-1}\bar
e_k+p\bar e_n.
$$

{\bf Cubes.} Any lattice cube is generated by some lattice point
$P$ and a $n$-tuple of linearly independent lattice vectors $\bar
v_i$:
$$
\Big\{ P+\sum\limits_{i=1}^n \alpha_i \bar v_i \Big| 0\le
\alpha_i\le 1, i=1,\ldots,n \Big\}
$$
We denote by $\{4,3^{n-2}\}^{L}_1$ for any $n\ge 2$ the lattice
cube with
a vertex at the origin and generated by all basis vectors.\\
By $\{4,3^{n-2}\}^{L}_2$ for any $n\ge 2$ we denote the lattice
cube with a vertex at the origin and generated by the first
$n{-}1$
basis vectors and the vector $\bar e_1+\bar e_2+\ldots+\bar e_{n-1}+2\bar e_{n}$.\\
By $\{4,3^{n-2}\}^{L}_3$ for any $n\ge 3$ we denote the lattice
cube with a vertex at the origin and generated by the vectors:
$\bar e_1$, and $\bar e_1 + 2\bar e_i$ for $i=2,\ldots,n$.

{\bf Generalized octahedra.} We denote by $\{3^{n-2},4\}^{L}_1$
for any $n\ge 2$ the lattice generalized octahedron
with the vertices $O\pm \bar e_i$ for $i=1,\ldots,n$.\\
By $\{3^{n-2},4\}^{L}_2$ for any positive $n$ we denote the
lattice generalized octahedron with the vertices $O\pm \bar e_i$
for $i=1,\ldots,n{-}1$,
and $O\pm\big(\bar e_1+\bar e_2+\ldots+\bar e_{n-1}+2\bar e_{n}\big)$.\\
By $\{3^{n-2},4\}^{L}_3$ for any positive $n$ we denote the
lattice generalized octahedron with the vertices $O$, $O-\bar
e_1$, $O-\bar e_1 - \bar e_i$ for $i=2,\ldots,n$, and $e_i$ for
$i=2,\ldots,n$.

{\bf A segment, octagons, and 24-sells.}
Denote by $\{\}^{L}$ the lattice segment with the vertices $O$ and $O+\bar e_1$.\\
By $\{6\}^{L}_1$ we denote the hexagon with the vertices $O\pm
\bar e_1$, $O\pm \bar e_2$,
$O\pm (\bar e_1-\bar e_2)$.\\
By $\{6\}^{L}_2$ we denote the hexagon with the vertices $O\pm
(2\bar e_1+\bar e_2)$,
$O\pm (\bar e_1+2\bar e_2)$, $O\pm (\bar e_1-\bar e_2)$.\\
By $\{3,4,3\}^{L}_1$ we denote the 24-sell with 8 vertices of the
form
$$
O\pm 2(\bar e_2+\bar e_3+\bar e_4), \quad O\pm 2(\bar e_1+\bar
e_2+\bar e_4), \quad O\pm 2(\bar e_1+\bar e_3+\bar e_4), \quad
O\pm 2\bar e_4,
$$
and 16 vertices of the form
$$
O\pm (\bar e_2+\bar e_3+\bar e_4)\pm (\bar e_1+\bar e_2+\bar e_4)
\pm (\bar e_1+\bar e_3+\bar e_4) \pm \bar e_4.
$$
By $\{3,4,3\}^{L}_2$ we denote the 24-sell with 8 vertices of the
form
$$
\begin{array}{l}
O\pm 2(\bar e_1+\bar e_2+\bar e_3+\bar e_4), \quad
O\pm 2(\bar e_1-\bar e_2+\bar e_3+\bar e_4), \\
O\pm 2(\bar e_1+\bar e_2-\bar e_3+\bar e_4), \quad O\pm 2(\bar
e_1+\bar e_2+\bar e_3-\bar e_4),
\end{array}
$$
and 16 vertices of the form
$$
O\pm (\bar e_1+\bar e_2+\bar e_3+\bar e_4) \pm (\bar e_1-\bar
e_2+\bar e_3+\bar e_4) \pm (\bar e_1+\bar e_2-\bar e_3+\bar e_4)
\pm (\bar e_1+\bar e_2+\bar e_3-\bar e_4).
$$

\subsection{Theorem on enumeration of convex elementary
lattice-regular lattice polytopes.} Now we formulate the main
statement of the work.

\begin{theorem}\label{main_theorem}
Any elementary lattice-regular convex lattice polytope is
lattice-con\-gru\-ent to some polytope of the following list.

{\bf List of the polygons.}\\
{\bf Dimension 1:} the segment $\{\}^{L}$.\\
{\bf Dimension 2:} the triangles $\{3\}^{L}_1$ and $\{3\}^{L}_2$;\\
the squares $\{4\}^{L}_1$ and $\{4\}^{L}_2$;\\
the octagons  $\{6\}^{L}_1$ and $\{6\}^{L}_2$.\\
{\bf Dimension 3:} the tetrahedra $\{3,3\}^{L}_i$, for $i=1,2,4$;\\
the octahedra $\{3,4\}^{L}_i$, for $i=1,2,3$;\\
the cubes $\{4,3\}^{L}_i$, for $i=1,2,3$.\\
{\bf Dimension 4:} the symplices $\{3,3,3\}^{L}_1$ and $\{3,3,3\}^{L}_5$;\\
the generalized octahedra $\{3,3,4\}^{L}_i$, for $i=1,2,3$;\\
the 24-sells $\{3,4,3\}^{L}_1$ and $\{3,4,3\}^{L}_2$;\\
the cubes $\{4,3,3\}^{L}_i$, for $i=1,2,3$.\\
{\bf Dimension n (n$>$4):} the symplices $\{3^{n-1}\}^{L}_i$
where positive integers
$i$ are divisors of~$n{+}1$;\\
the generalized octahedra $\{3^{n-2},4\}^{L}_i$, for $i=1,2,3$;\\
the cubes $\{4,3^{n-2}\}^{L}_i$, for $i=1,2,3$.

All polytopes of this list are lattice-regular. Any two polytopes
of the list are not lattice-congruent to each other.
\end{theorem}

\begin{figure}[h]
$$\epsfbox{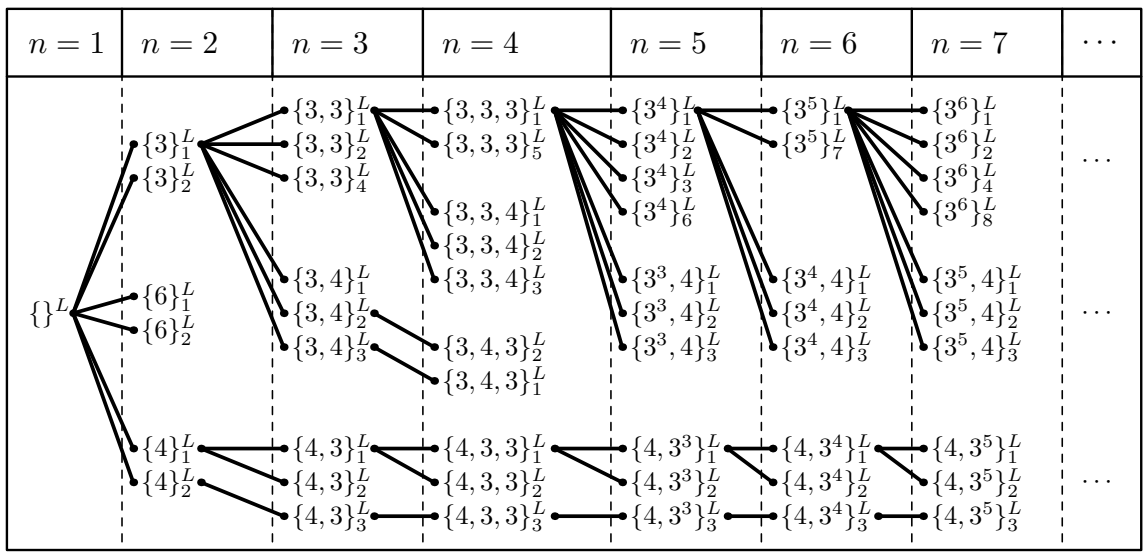}$$
\caption{The adjacency diagram for the elementary lattice-regular
convex lattice polytopes.}\label{adjac}
\end{figure}

On Figure~\ref{adjac} we show the adjacency diagram for the
elementary lattice-regular convex lattice polygons of dimension
not exceeding 7. Lattice-regular lattice polygons of different
(six) types are shown on Figure~\ref{pic_reg} in the next section.
Lattice-regular lattice three-dimensional polygons of different
(nine) types are shown on Figures~\ref{3dsymp}, \ref{3dcube},
and~\ref{3doct} further in Sections~\ref{symplex}, \ref{cube},
and~\ref{octahedron} respectively.

\vspace{2mm}

Further in the proofs we will use the following definition.
Consider a $k$-dimensional lattice polytope $P$. Let its
Euclidean volume equal $V$. Denote the Euclidean volume of the
minimal $k$-dimensional symplex in the $k$-dimensional plane of
the polytope by $V_0$. The ratio $V/V_0$ is said to be the {\it
lattice volume} of the given polytope (if $k=1$, or $2$ --- the
{\it lattice length} of the segment, or the {\it lattice area} of
the polygon respectively).

\section{Two-dimensional case.}\label{2d}

In this section we prove Theorem~\ref{main_theorem} for the
two-dimensional case.

\begin{proposition}\label{regular_poligons}
Any elementary lattice-regular $($two-dimensional$)$ lattice
convex polygon
is lattice-congruent to one of the following polygons $($see Figure~\ref{pic_reg}$)$:\\
{\bf 1)} $\{3\}^{\z}_1$; {\bf 2)} $\{3\}^{\z}_2$; {\bf 3)}
$\{4\}^{\z}_1$; {\bf 4)} $\{4\}^{\z}_2$; {\bf 5)} $\{6\}^{\z}_1$;
{\bf 6)} $\{6\}^{\z}_2$.
\begin{figure}[h]
$$\epsfbox{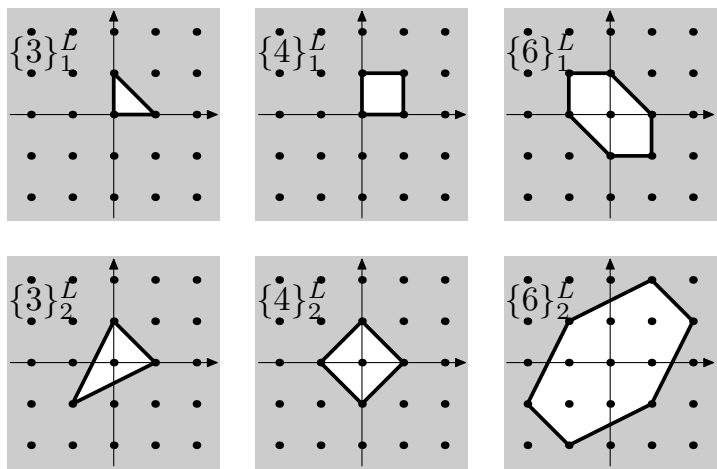}$$
\caption{The lattice-regular polygons with edges of unit
length.}\label{pic_reg}
\end{figure}
\end{proposition}

\begin{proof}
Suppose that the lattice polygon $A_1A_2\ldots A_n$, where $n\ge
3$, is primitive and lattice-regular.

Let us prove that all edges of $A_1A_2\ldots A_n$ are of unite
lattice lengths. Since the polygon is lattice-regular, all its
edges are lattice-congruent, and hence they are of the same
lattice length. Suppose that the lattice lengths of all edges
equal some positive integer $k$. Then our polygon is
lattice-congruent to the $k$-tuple of the polygon $A'_1A'_2\ldots
A'_n$, where $A'_1=A_1$, and $A'_l=A_l{+}1/k(\bar{A_{l}A_{l+1}})$
for $l=2, \ldots, n{-}1$. Therefore, $k=1$.

Denote by $B_i$ the midpoint of the edge $A_iA_{i+1}$ for
$i=1,\ldots ,n{-}1$, and by $B_n$ the midpoint of the edge
$A_nA_1$. Suppose that $n$ is even ($n=2\tilde n$). Denote by $M$
the midpoint of the segment $A_1A_{\tilde n+1}$. Note that the
point $M$ is the common intersection of the segments
$A_iA_{\tilde n+i}$ for $i=1,\ldots, \tilde n$. Suppose that $n$
is odd ($n=2\tilde n+1$). Denote by $M$ the intersection point of
the segments $A_1B_{\tilde n}$ and $A_2B_{\tilde n+1}$. Note that
the point $M$ is the common intersection of the segments
$A_iB_{\tilde n+i}$, $A_{\tilde n+i+1}B_i$ for $i=1,\ldots,
\tilde n$, and the segment $A_{\tilde n}B_{2\tilde n+1}$.

For any integer $i$ such that $1\le i\le n$ the following holds.
The transformation that preserves the points $M$ and $B_i$, and
taking the point $A_i$ to the point $A_{i+1}$ (or $A_n$ to $A_1$
in the case of $i=n$) is lattice-affine and preserve the polygon
$A_1A_2\ldots A_n$.

Suppose that the polygon $A_1A_2\ldots A_n$ contains some lattice
point not contained in the union of its vertices and segments
$MB_i$ for $i=1,\ldots, n$. Then by symmetry reasons the triangle
$ A_1MB_1$ contains at least one lattice point, that is not
contained in the edges $A_1B_1$ and $MB_1$. Denote one of such
points by $P$. Let $Q$ be the point symmetric to the point $P$
about the line $MB_1$. The segment $PQ$ is parallel to the
segment, and hence the lattice point $A_1 + \bar {PQ}$ is
contained in the interior of the segment $A_1A_2$. Then the
lattice length of the edge $A_1A_2$ is not unit. We come to the
contradiction with the above.

Therefore, all inner lattice points of the polygon $A_1A_2\ldots
A_n$ are contained in the union of the segments $MB_i$ for
$i=1,\ldots, n$ and vertices. Now we study all different cases of
configurations of lattice points on the segment $MB_1$.

\vspace{2mm}

{\it Case 1.} Suppose that $MB_1$ does not contain lattice points.
Then by symmetry reasons the polygon $A_1A_2\ldots A_n$ does not
contain lattice points different to its vertices. Hence the
vectors $\bar{A_2A_1}$ and $\bar{A_2A_3}$ generate the lattice.
Consider the linear system of coordinates such that the points
$A_1$, $A_2$, and $A_3$ have the coordinates $(0,1)$, $(0,0)$,
and $(1,0)$ in it respectively.

If $n=3$, then the triangle $A_1A_2A_3$ is lattice-congruent to
the triangle $\{3\}^L_1$.

Let $n>3$. Since the vectors $\bar {A_1A_2}$ and $\bar {A_2A_3}$
generate the lattice, and the vectors $\bar {A_2A_3}$ and $\bar
{A_3A_4}$ generate the lattice, the point $A_4$ has the
coordinates $(a,1)$ for some integer $a$. Since the segment
$A_1A_4$ does not contain lattice points distinct to the
endpoints, $A_4=(1,1)$. By the same reasons $A_n=(1,1)$.
Therefore, $n=4$, and the lattice polygon $A_1A_2A_3A_4$ is
lattice-congruent to the lattice-regular quadrangle $\{4\}^L_1$.

\vspace{2mm}

{\it Case 2.} Suppose that the point $M$ is lattice and the
segment $MB_1$ does not contain lattice points distinct to $M$.
Then the vectors $\bar{MA_1}$ and $\bar{MA_2}$ generate the
lattice. Consider the linear system of coordinates such that the
points $A_1$, $M$, and $A_2$ have the coordinates $(0,1)$,
$(0,0)$, and $(1,0)$ in it. Since the vectors $\bar {A_1M}$ and
$\bar {MA_2}$ generate the lattice, and the vectors $\bar {A_2M}$
and $\bar {MA_3}$ generate the lattice, the point $A_3$ has the
coordinates $(-1,a)$ for some integer $a$.

If $a\ge 2$, then the polygon is not convex or it contains
straight angles.

If $a=1$, then the vectors $\bar {A_1A_2}$ and $\bar {A_2A_3}$
generate the lattice. Since the vectors $\bar {A_3M}$ and $\bar
{MA_4}$ generate the lattice, and the vectors $\bar {A_2A_3}$ and
$\bar {A_3A_4}$ generate the lattice, the new coordinates of the
point $A_4$ are $(0,-1)$. Since $A_4=A_1{+}2\bar{A_1M}$, we have
$$
A_5=A_2+2\bar{A_2M}=(0,-1), \qquad A_6=A_3+2\bar{A_3M}=(1,-1),
\quad \hbox{and} \quad n=6.
$$
Therefore, the lattice-regular polygon $A_1A_2A_3A_4A_5A_6$ is
lattice-congruent to the lattice-regular hexagon $\{6\}^L_1$.

If $a=2$, then $A_3=A_1{+}2\bar{A_1M}$. Hence
$A_3=A_2{+}2\bar{A_2M}=(0,-1)$, and $n=4$. Therefore, the polygon
$A_1A_2A_3A_4$ is lattice-congruent to the lattice-regular
quadrangle $\{4\}^L_2$.

If $a=3$, then $A_3$ is contained in the line $MB_1$. Hence $n=3$.
Therefore, the lattice triangle $\triangle A_1A_2A_3$ is
lattice-congruent to the lattice-regular triangle $\{3\}_2$.

Since $A_n=(a,-1)$, the edges $A_nA_1$ and $A_1A_2$ intersect for
the case of $a>3$.

{\it Case 3.} Suppose that the segment $MB_1$ contains the unique
lattice point $P$ distinct from the endpoints of the segment
$MB_1$. Then the vectors $\bar{PA_1}$ and $\bar{PA_2}$ generate
the lattice. Consider the linear system of coordinates such that
the points $A_1$, $P$, and $A_2$ have the coordinates $(0,1)$,
$(0,0)$, and $(1,0)$ in it. Since the polygon is lattice-regular,
the point $M{+}\bar{PM}$ is also a lattice point, and hence
$M=(-1/2,-1/2)$. Denote the point (-1,-2) by $M'$. (Note  that
the point $M$ is the midpoint of the segment $M'B$).

The vectors $\bar {A_1M'}$ and $\bar {M'A_2}$ generate a
sublattice of index $3$. The vectors $\bar {A_2M'}$ and $\bar
{M'A_3}$ generate a sublattice of index $3$.  The segment
$A_2A_3$ is of unit lattice length. Therefore, the point $A_3$
has the coordinates $(2a{-}1,6a{+}2)$ for some integer $a$.

If $a\ge 0$, then the polygon is not convex, or it contains
straight angles. Since $A_n=(6a{+}2,2a{-}1)$, the edges $A_nA_1$
and $A_1A_2$ intersect for the case of $a<0$.

{\it Case 4.} Suppose that the point $M$ is lattice and the
segment $MB_1$ contains the unique interior lattice point $P$.
Then the vectors $\bar{PA_1}$ and $\bar{PA_2}$ generate the
lattice. Consider the linear system of coordinates such that the
points $A_1$, $P$, and $A_2$ have the coordinates $(0,1)$,
$(0,0)$, and $(1,0)$ in it. Since the polygon is lattice-regular,
the point $M+\bar{PM}$ is also lattice, and hence $M=(-1,-1)$.

The vectors $\bar {A_1M'}$ and $\bar {M'A_2}$ generate a
sublattice of index $3$. The vectors $\bar {A_2M'}$ and $\bar
{M'A_3}$ generate a sublattice of index $3$.  The segment
$A_2A_3$ is of unit lattice length. Therefore, the point $A_3$
has the coordinates $(a{-}1,2a{+}2)$ for some integer $a$, such
that $a\not\equiv 1(\MOD 3)$.

If $a\ge 0$, then the polygon is not convex, or it contains
straight angles, but this is impossible.

If $a=-1$, then the vectors $\bar {A_1A_2}$ and $\bar {A_2A_3}$
generate a sublattice of index $3$. Since the vectors $\bar
{A_3M'}$ and $\bar {M'A_4}$ generate a sublattice of index $3$,
and the vectors $\bar {A_2A_3}$ and $\bar {A_3A_4}$ generate a
sublattice of index $3$, the point $A_4=(-3,-2)$. Since
$A_4=A_1{+}2\bar{A_1M}$, we have
$$
A_5=A_2+2\bar{A_2M}=(-2,-3), \qquad A_6=A_3+2\bar{A_3M}=(0,-2),
\quad \hbox{and} \quad n=6.
$$
Therefore, the lattice polygon $A_1A_2A_3A_4A_5A_6$ is
lattice-congruent to the lattice-regular hexagon $\{6\}^L_2$.

Since $A_n=(2a{+}2,a{-}1)$, the edges $A_nA_1$ and $A_1A_2$ are
intersecting for the case of $a<-1$.

{\it The remaining cases.} Suppose that the segment $MB_1$
contains at least two interior lattice points. Let $P_1$ and
$P_2$ be two distinct interior lattice points of the segment
$MB_1$. Let also the segment $MP_2$ contains the point $P_1$.

Consider an lattice-affine transformation $\xi$ taking the point
$M$ to itself, and the segment $A_1A_2$ to the segment $A_2A_3$.
The points $Q_1=\xi(P_1)$ and $Q_2=\xi(P_2)$ are contained in the
segment $MB_2$. Since the lines $P_1Q_1$ and $P_2Q_2$ are parallel
and the triangle $P_2MQ_2$ contains the segment $P_1Q_1$, the
lattice point $S=P_2{+}\bar{Q_2P_2}$ is contained in the interior
of the segment $P_1Q_1$. Hence the lattice point $S$ of the
polygon $A_1A_2\ldots A_n$ is not contained in the union of
segments $MB_i$ for $i=1,\ldots, n$. We come to the contradiction.

\vspace{2mm}

We have studied all possible cases of configurations of lattice
points contained inside lattice polygons. The proof of
Proposition~\ref{regular_poligons} is completed.
\end{proof}

\section{Lattice-regular lattice symplices.}\label{symplex}

In this section we study all lattice-regular lattice symplices
for all integer dimensions.

Let us fix some basis of lattice vectors $\bar e_i$, for
$i=1,\ldots,n$ generating the lattice in $\r^n$ and the
corresponding coordinate system. Denote by $O$ the origin in
$\r^n$.

\begin{proposition}\label{regular_symplices}
{\bf Sym$_1$.} All elementary lattice-regular one-dimensional
lattice symplices are lattice segments of unit lattice length.

{\bf Sym$_n$ (for $n>1$).} {\it i$)$} The symplex
$\{3^{n-1}\}^L_n$ where $p$ is a positive divisor of $n{+}1$ is
elementary and lattice-regular;
\\
{\it ii$)$} any two symplices listed in {\it $($i$)$} are not
lattice-congruent to each other;
\\
{\it iii$)$} any elementary lattice-regular $n$-dimensional
lattice symplex is lattice-congruent to one of the symplices
listed in {\it $($i$)$}.
\end{proposition}

\begin{figure}[h]
$$\epsfbox{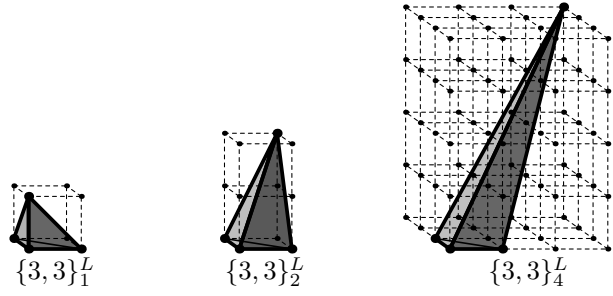}$$
\caption{Three-dimensional elementary lattice-regular convex
lattice tetrahedra.}\label{3dsymp}
\end{figure}

The three-dimensional tetrahedra $\{3,3\}^L_1$, $\{3,3\}^L_2$,
and $\{3,3\}^L_4$ are shown on Figure~\ref{3dsymp}.

\begin{proof}
We start the proof of studying of some low-dimensional cases. The
one-dimensional case is trivial and is omitted here. The
two-dimensional case was described in
Proposition~\ref{regular_poligons}. Let us study the
three-dimensional case.

\vspace{2mm}

{\it Three-dimensional case.} Consider an arbitrary elementary
lattice-regular three-di\-men\-si\-onal lattice tetrahedron $S$.
Since its faces are lattice-regular, by
Proposition~\ref{regular_poligons} the faces are
lattice-congruent either to $\{3\}^L_1$ or to $\{3\}^L_3$.

Suppose that the faces of $S$ are lattice-congruent to
$\{3\}^L_1$. Then there exist a positive integer $b$, nonnegative
integers $a_1$, $a_2$ less than $b$, and a lattice-affine
transformation taking the the tetrahedron $S$ to the tetrahedron
$S'$ with the vertices
$$
V_0=O, \quad V_1=O+\bar e_1, \quad V_2=O+\bar e_2, \quad
V_3=O+a_1\bar e_1+a_2\bar e_2+b\bar e_3.
$$

Since $S'$ is also a lattice-regular tetrahedron, the group of
its symmetries is isomorphic to the group of permutations of
order 4. This group is generated by the following transpositions
of vertices: $V_1$ and $V_2$, $V_2$ and $V_3$, and $V_0$ and
$V_2$. The first two transpositions are linear, and the third one
is linear after shifting by the vector $-\bar e_1$. Direct
calculations shows, that the matrices of the corresponding linear
transformations are the following:
$$
\left(
\begin{array}{ccc}
0&1&\frac{a_1-a_2}{b}\\
1&0&\frac{a_2-a_1}{b}\\
0&0&1\\
\end{array}
\right), \quad \left(
\begin{array}{ccc}
1&a_1&-\frac{a_1(a_2+1)}{b}\\
0&a_2&\frac{1-a_2^2}{b}\\
0&b&-a_2\\
\end{array}
\right), \quad \hbox{and} \quad \left(
\begin{array}{ccc}
1&0&0\\
-1&-1&\frac{2a_2+a_1-1}{b}\\
0&0&1\\
\end{array}
\right).
$$
Since the listed transformations are lattice-linear we have only
the following possibilities, all these matrices are integer.
Therefore, $a_1\equiv a_2 \equiv b{-}1(\MOD b)$. Since, positives
$a_1$ and $a_2$ was chosen to be smaller than $b$, we have
$a_1=a_2=b{-}1$. Since the matrix $A_3$ is integer, the
coefficient $3{-}4/b$ is also integer. So we have to check only
the following cases for $a_1$, $a_2$, and $b$: $b=1$, and
$a_1=a_2=0$; $b=2$, and $a_1=a_2=1$; $b=4$, and $a_1=a_2=3$.
These cases corresponds to the tetrahedra $\{3,3\}^L_1$,
$\{3,3\}^L_2$, and $\{3,3\}^L_4$ respectively. Since the lattice
volume of $\{3,3\}^L_p$ equals $p$, the above tetrahedra are not
lattice-congruent.

Let us prove that the faces of $S$ are not lattice-congruent to
$\{3\}^L_3$ by reductio ad absurdum. Suppose it is so. Let $V_0$,
$V_1$, $V_2$, and $V_3$ be the vertices of $S$. Since the faces
$V_0V_2V_3$ and $V_1V_2V_3$ are congruent to $\{3\}^L_3$, the
face $V_0V_2V_3$ contains a unique lattice point in its interior
(we denote it by $P_1$), and the face $V_1V_2V_3$ contains a
unique lattice point in its interior (we denote it by $P_2$).
Consider a lattice-symmetry of $S$ permuting $V_0$ and $V_1$ and
preserving $V_2$ and $V_3$. This symmetry takes the face
$V_0V_2V_3$ to the face $V_1V_2V_3$, and hence it maps the point
$P_1$ to $P_2$. Thus, the lattice vector $\bar{P_1P_2}$ is
parallel to the vector $V_0V_1$. Hence, the point
$V_0+\bar{P_1P_2}$ is interior lattice point of the segment
$V_0V_1$. Therefore, the segment $V_0V_1$ is not of unit lattice
length and the face $V_0V_1V_2$ is not lattice-congruent to
$\{3\}^L_3$. We come to the contradiction.

This completes the proof of Proposition~\ref{regular_symplices}
for the three-dimensional case.

\vspace{2mm}

In the higher dimensional case we first study two certain
families of symplices.

\vspace{2mm}

{\it The first type of $n$-dimensional symplices.} Consider any
symplex (for $n>3$) with the vertices
$$
V_0=O, \qquad \hbox{$V_k=O+\bar e_i$, for $k=1,\ldots,n{-}1$},
\quad \hbox{and} \quad V_n=O+\sum\limits_{i=1}^n a_i\bar e_i,
$$
we denote it by $S^n(a_1,\ldots, a_n)$. We suppose that all $a_i$
are nonnegative integers satisfying $a_k<a_n$ for $k=1,\ldots,
n{-}1$. Let us find the conditions on $a_i$ for the symplex to be
lattice-regular.

If $S^n(p;a_1,\ldots, a_n)$ is a lattice-regular symplex then the
group of its symmetries is isomorphic to the group of
permutations of order $n{+}1$. This group is generated by the
following transpositions of vertices: the transposition
exchanging $V_k$ and $V_{k+1}$ for $k=1,\ldots, n{-}1$, and the
transposition exchanging $V_0$ and $V_2$. The first $n{-}1$
transpositions are linear (let their matrices be $A_k$ for
$k=1,\ldots, n{-}1$), and the last one is linear after shifting
by the vector $-\bar e_1$ (denote the corresponding matrix  by
$A_n$). Let us describe the matrices of these transformations
explicitly.

The matrix $A_k$ for $k=1, \ldots n{-}2$ coincides with the
matrix transposing the vectors $\bar e_k$ and $\bar e_{k+1}$,
except the last column. The $n$-th column contains the
coordinates of the vector
$$
\begin{array}{c}
\frac{a_k-a_{k+1}}{a_n}(\bar e_k-\bar e_{k+1}) +\bar e_n.\\
\end{array}
$$
The matrix $A_{n-1}$ coincides with the matrix of identity
transformation except the last two columns. These columns contain
the coefficients of the following two vectors respectively:
$$
\begin{array}{c}
\sum\limits_{j=1}^{n}a_j\bar e_j, \quad \hbox{and} \quad
\sum\limits_{j=1}^{n-2}\left( \frac{-a_j(a_{n-1}+1)}{a_n}\bar
e_j\right) + \frac{1-a_{n-1}^2}{a_n}\bar e_{n-1}-a_{n-1}\bar e_n.
\end{array}
$$
The matrix $A_{n}$ coincides with the unit matrix except the
second row. This row is as follows: $ (-1,\ldots, -1,
\frac{a_1+\ldots+a_{n-1}-1+a_2}{a_n}). $

The determinants of all such matrices equal $-1$. So the
corresponding affine transformations are lattice iff all the
coefficients of all the matrices are lattice. The matrices $A_k$
for $k\le n{-}2$ are lattice iff
$$
a_1\equiv a_2\ldots \equiv a_{n-1} (\MOD a_n).
$$
Since $a_k<a_n$, we have the equalities. Suppose
$a_1=\ldots=a_{n-1}=p{-}1$ for some positive integer $p$. The
matrix $A_{n-1}$ is integer, iff $1{-}p^2\equiv p(p{+}1)\equiv 0
(\MOD a_n)$. Therefore, $r{+}1$ is divisible by $a_n$, and hence
$a_n=p$. The matrix $A_n$ is integer, iff $n(p{-}1){-}1$ is
divisible by $p$, or equivalently $n{+}1$ is divisible by $p$.

So, we have already obtained that the symplex $S^n(a_1,\ldots,
a_n)$ where $n>3$ and $0\le a_k<a_n$ for $k=1,\ldots, n{-}1$ is
lattice-regular iff it is coincides with some $\{3^{n-1}\}^L_p$
for some $p$ dividing $n$. Since the lattice volume of $S_p^n$
equals $p$, the above symplices are not lattice-congruent.

\vspace{2mm}

{\it The second type of $n$-dimensional symplices.} Here we study
symplices (for $n>3$) with the vertices
$$
\begin{array}{c}
V_0=O, \qquad \hbox{$V_k=O+\bar e_i$, for $k=1,\ldots,n{-}2$},\\
\qquad V_{n-1}=O+(p-1)\sum\limits_{i=1}^{n-2}\bar e_i + p\bar
e_{n-1}, \quad \hbox{and} \quad V_n=O+\sum\limits_{i=1}^n a_i\bar
e_i,
\end{array}
$$
denote such symplices by $S^n(p;a_1,\ldots, a_n)$. We also
suppose, that all $a_i$ are nonnegative integers satisfying
$a_k<a_n$ for $k=1,\ldots, n{-}1$, and $p\ge 2$. Let us show that
all these symplices are not lattice-regular. Consider an
arbitrary symplex $S^n(p;a_1,\ldots, a_n)$, satisfying the above
conditions.

Consider the symmetry exchanging $V_{n-1}$ and $V_{n}$. This
transformation is linear. Its matrix coincide with the matrix of
identity transformation, except for the last two columns. These
columns contains the coefficients of the following two vectors
respectively:
$$
\begin{array}{l}
\sum\limits_{j=1}^{n-2}\left(\frac{a_j+1}{p}-1\right)\bar e_j+
\frac{a_{n-1}}{p}\bar e_{n-1}+\frac{a_{n}}{p}\bar e_{n},
\quad \hbox{and}\\
\sum\limits_{j=1}^{n-2}\left(
\frac{p(a_{n-1}-a_j+p-1)-a_{n-1}-a_ja_{n-1}}{pa_n}\bar e_j\right)
+ \frac{p^2-a_{n-1}^2}{pa_n}\bar e_{n-1}-\frac{a_{n-1}}{p}\bar
e_n.
\end{array}
$$
If this transformation is lattice-linear, then $a_2{+}1$ is
divisible by $p$.

Consider the symmetry exchanging $V_{1}$ and $V_{n}$. This
transformation is linear. Its matrix coincide with the matrix of
identity transformation, except for the first column and the last
two columns. These columns contains the coefficients of the
following three vectors respectively:
$$
\begin{array}{c}
\sum\limits_{j=1}^{n}a_j\bar e_j, \qquad
-\sum\limits_{j=1}^{n}\left(\frac{a_j(p-1)}{p}\bar e_j\right)+
\frac{p-1}{p}\bar e_1+\bar e_{n-1},
\quad \hbox{and}\\
\frac{a_{n-1}p-pa_1-a_{n-1}-p}{pa_n}\sum\limits_{j=1}^{n}a_j\bar
e_j+ \frac{p+a_1}{pa_n}\bar e_1+\bar e_n.
\end{array}
$$
If this transformation is lattice-linear then $a_2$ is divisible
by $r$.

Since $a_2$ and $a_2{+}1$ are divisible by $r$ and $r\ge 2$, the
symplex $S^n(p;a_1,\ldots, a_n)$ is not lattice-regular.

\vspace{2mm}

{\it Conclusion of the proof.}

Statement {\it $($i$)$} is already proven. Since $p$ is a lattice
volume of $S_p^n$, Statement {\it $($ii$)$} holds. We prove
Statement {\it $($iii$)$} of the proposition by the induction on
the dimension $n$. For $n=1,2,3$ the statement is already proven.
Suppose that it is true for an arbitrary $n \ge 3$. Let us prove
the statement for $n{+}1$.

Consider any lattice-regular $(n{+}1)$-dimensional lattice
symplex $S$. Since it is lattice-regular, all its faces are
lattice-regular. By the induction assumption there exist a
positive integer $p$ dividing $n{+}1$ such that the faces of $S$
are lattice-congruent to $\{3^{n-1}\}^L_p$. Therefore, $S$ is
lattice-affine equivalent to the symplex $S^{n+1}(p;a_1,\ldots,
a_{n+1})$, where $p\ge 1$, and $a_k<a_n$ for $k=1,\ldots, n{-}1$.
By the above cases the lattice-regularity implies, that $p=1$,
and that there exist a positive integer $p'$ dividing $n{+}2$
such that the symplex $S^{n+1}(1;a_1,\ldots, a_{n+1})$ coincides
with $\{3^{n}\}^L_{p'}$. This concludes the proof of the
Statement {\it $($iii$)$} for the arbitrary dimension.

Proposition~\ref{regular_symplices} is proven.
\end{proof}

\section{Lattice-regular lattice cubes.}\label{cube}

In this section we describe all lattice-regular lattice cubes for
all integer dimensions.

\begin{proposition}\label{regular_cubes}
{\bf Cube$_1$.} All elementary lattice-regular one-dimensional
lattice cubes are lattice segments of unit lattice length.

{\bf Cube$_2$.} All elementary lattice-regular two-dimensional
lattice cubes are lattice-con\-gru\-ent to $\{4\}^L_1$, or to
$\{4\}^L_2$. The cubes $\{4\}^L_1$ and $\{4\}^L_2$ are not
lattice-congruent.

{\bf Cube$_n$ (for $n>2$).} All elementary lattice-regular
$n$-dimensional lattice cubes are lat\-tice-con\-gru\-ent to
$\{4,3^{n{-}2}\}^L_1$, $\{4,3^{n{-}2}\}^L_2$, or
$\{4,3^{n{-}2}\}^L_3$. The cubes $\{4,3^{n{-}2}\}^L_1$,
$\{4,3^{n{-}2}\}^L_2$, and $\{4,3^{n{-}2}\}^L_3$ are not
lattice-congruent to each other.
\end{proposition}

\begin{figure}[h]
$$\epsfbox{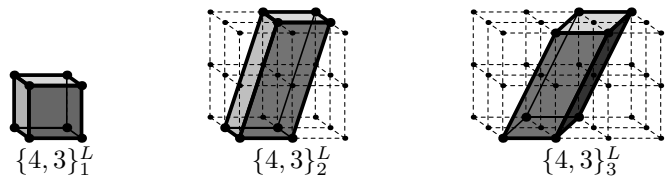}$$
\caption{Three-dimensional elementary lattice-regular convex
lattice cubes.}\label{3dcube}
\end{figure}

The three-dimensional cubes $\{4,3^{n{-}2}\}^L_1$,
$\{4,3^{n{-}2}\}^L_2$, and $\{4,3^{n{-}2}\}^L_3$ are shown on
Figure~\ref{3dcube}.

We use the following two facts.

\begin{lemma}\label{interior_points_cube}
Any elementary lattice-regular $n$-dimensional lattice cube
contains at most one lattice point in its interior. If the cube
contains an interior lattice point, this point coincides with the
intersection point of the diagonals of the cube.
\end{lemma}

\begin{proof}
Suppose the converse is true. There exist a primitive
$n$-dimensional cube with an interior lattice point $A$ distinct
to the intersections of the diagonals. Then there exist a
lattice-linear reflection $\xi$ of the cube above some $n{-}1$
dimensional plane, that do not preserve $A$. So the line of
vector $\bar {A\xi(A)}$ coincide with the line containing one of
the generating vectors of the cube. Hence, the one-dimensional
faces of the cube are not elementary. Therefore, the cube is not
elementary. We come to the contradiction.
\end{proof}

Let us give the following important definition. Consider some
hyperplane containing a sublattice of the lattice of corank 1 and
a lattice point in the complement to this hyperplane. Let the
Euclidean distance from the given point to the given hyperplane
equals $l$. The minimal value of nonzero Euclidean distances from
the points of the lattice to the hyperplane is denoted by $l_0$.
The ratio $l/l_0$ is said to be the {\it lattice distance} from
the given lattice point to the given lattice hyperplane.

\begin{corollary}\label{ind_dis_cube}
The lattice distances from any vertex of any elementary
lattice-regular $n$-dimensional lattice cube to any its
$(n{-}1)$-dimensional face $($that does not containing the given
vertex$)$ equals either $1$, or $2$. \qed
\end{corollary}

{\it Proof of Proposition~\ref{regular_cubes}.} The
one-dimensional case is trivial. The two-dimensional case was
described in Proposition~\ref{regular_poligons}. Let us study
higher-dimensional cases.

\vspace{2mm}

{\it The lattice cubes $\{4,3^{n-2}\}^L_1$, $\{4,3^{n-2}\}^L_2$,
and $\{4,3^{n-2}\}^L_3$.} First, let us study the cases of the
polytopes $\{4,3^{n-2}\}^L_1$, $\{4,3^{n-2}\}^L_2$, and
$\{4,3^{n-2}\}^L_3$ for any $n\ge 3$. Since lattice volumes of
$\{4,3^{n-2}\}^L_1$, $\{4,3^{n-2}\}^L_2$, and $\{4,3^{n-2}\}^L_3$
are $n!$, $2n!$, and $2^{n-1}n!$ respectively, the listed cubes
are not lattice-congruent to each other. Let us prove that these
polytopes are lattice-regular for any $n\ge 3$.

Since the vectors of $\{4,3^{n-2}\}^L_1$ generate lattice, it is
lattice-regular.

Now we study the case of $\{4,3^{n-2}\}^L_2$. Denote by $\bar
v_i$ the vector $\bar e_i$ for $i=1,\ldots, n{-}1$ and by $\bar
v_n$ the vector $\bar e_1 +\bar e_2+ \ldots +\bar e_{n-1} + 2\bar
e_n$. The group of lattice symmetries of the cube
$\{4,3^{n-2}\}^L_2$ is generated by the linear operators $A_k$
transposing the vectors $\bar v_i$ and $\bar v_{i+1}$ for
$i=1,\ldots, n{-}1$, and the last one: the composition of the
symmetry $A_n$ sending $\bar v_1$ to $-\bar v_1$ and preserving
$\bar v_i$ for $i=2,\ldots, n$ and the lattice shift on the
vector $\bar v_1$.

Let us check that all the linear transformations $A_k$ are
lattice-linear. We show explicitly the matrices of $A_k$ in the
basis $\bar e_i$ for $i=1,\ldots, n$. The matrix of $A_k$ for
$k=1,\ldots, n{-}2$ coincides with the the matrix of the
transposition of the vectors $\bar e_k$ and $\bar e_{k+1}$. The
matrix of $A_{n-1}$ coincides with the matrix of identity
transformation except the last two columns. These columns contain
the coefficients of the vectors $\bar v_n$ and $\bar e_{n-1}+\bar
e_n-\bar v_n$ respectively. The matrix of $A_n$ coincides with
the matrix of identity transformation except the first row, which
is $(-1,0,\ldots,0,1)$. Since all these matrices are in
$SL(n,\z)$, the cube $\{4,3^{n-2}\}^L_2$ is lattice-regular.

\vspace{2mm}

Let us consider now the case of $\{4,3^{n-2}\}^L_3$. Put by
definition $\bar v_1=\bar e_1$, and $\bar v_i=\bar e_1 + 2\bar
e_i$
 for $i=2,\ldots,n$.
The group of lattice symmetries of the cube $\{4,3^{n-2}\}^L_3$
is generated by the linear operators $A_k$ transposing the vectors
$\bar v_i$ and $\bar v_{i+1}$ for $i=1,\ldots, n{-}1$, and the
last one: the composition of the symmetry $A_n$ sending $\bar
v_1$ to $-\bar v_1$ and preserving $\bar v_i$ for $i=2,\ldots,
n$, and the lattice shift on the vector $\bar v_1$.

Let us check that all the linear transformations $A_k$ are
lattice-linear. The matrix of $A_1$ coincides with the matrix of
identity transformation except the second row, which is
$(2,-1,\ldots, -1)$. The matrix of $A_k$ for $k=2,\ldots, n{-}1$
coincides with the matrix transposing the vectors $\bar e_k$ and
$\bar e_{k+1}$. The matrix of $A_n$ coincides with the matrix of
identity transformation except the first row, which is
$(-1,1,\ldots, 1)$. Since all these matrices are in $SL(n,\z)$,
the cube $\{4,3^{n-2}\}^L_3$ is lattice-regular.

\vspace{2mm}

{\it Conclusion of the proof of Proposition~\ref{regular_cubes}.}
Now we prove, that any elementary lattice-regular $n$-dimensional
lattice cube for $n \ge 3$ is lattice-congruent to
$\{4,3^{n-2}\}^L_1$, $\{4,3^{n-2}\}^L_2$, or $\{4,3^{n-2}\}^L_3$.
by the induction on $n$.

{\it The base of induction.} Any face of any three-dimensional
lattice-regular cube is lattice-regular.

Suppose, that the faces of three-dimensional lattice-regular cube
$C$ are lattice-congruent to $\{4\}^L_1$. Then $C$ is
lattice-congruent to the cube generated by the origin and the
vectors $\bar e_1$, $\bar e_2$, and $a_1\bar e_1+a_2\bar
e_2+a_3\bar e_3$. By Corollary~\ref{ind_dis_cube} we can choose
$a_3$ equals either $1$ or $2$. If $a_3=1$ then $C$ is
lattice-congruent to $\{4,3\}^L_1$. If $a_3=2$, then we can
choose $a_1$ and $a_2$ being $0$, or $1$. Direct calculations
show, that the only possible case is $a_1=1$, and $a_2=0$, i.e.
the case of $\{4,3\}^L_2$.

Suppose now,  that the faces of three-dimensional lattice-regular
cube $C$ are lattice-congruent to  $\{4\}^L_2$. Then $C$ is
lattice-congruent to the cube generated by the origin and the
vectors $\bar e_1$, $\bar e_1+2\bar e_2$, and $a_1\bar
e_1+a_2\bar e_2+a_3\bar e_3$. By Corollary~\ref{ind_dis_cube} we
can choose $a_3$ equals either $1$ or $2$. Then we can also
choose $a_1$ and $a_2$ being $0$, or $1$. Direct calculations
show, that the only possible case $a_1=1$, $a_2=0$, and $a_3=2$
corresponds to $\{4,3\}^L_3$.

{\it The step of induction.} Suppose that any elementary
lattice-regular $(n{-}1)$-dimensional lattice cubes ($n > 3$) are
lattice-congruent to $\{4,3^{n-3}\}^L_1$, $\{4,3^{n-3}\}^L_2$, or
$\{4,3^{n-3}\}^L_3$. Let us prove that any elementary
lattice-regular $n$-dimensional lattice cubes are
lattice-congruent to $\{4,3^{n-2}\}^L_1$, $\{4,3^{n-2}\}^L_2$, or
$\{4,3^{n-2}\}^L_3$.

Any face of any lattice-regular cube is lattice-regular. Suppose,
that the faces of $(n{-}1)$-dimensional lattice-regular cube $C$
are lattice-congruent to $\{4,3^{n-3}\}^L_1$. Then $C$ is
lattice-congruent to the cube $C'$ generated by the origin and
the vectors $\bar v_i =\bar e_i$ for $i=1,\ldots, n-1$ and the
vector $\bar e_n=a_1\bar e_1+\ldots+a_n\bar e_n$.

By Corollary~\ref{ind_dis_cube} we can choose $a_n$ equals either
$1$ or $2$. If $a_n=1$, then the lattice volume of $C$ is $n!$ and
it is lattice-congruent to $\{4,3^{n-2}\}^L_1$. If $a_n=2$, then
we can choose $a_i$ being $0$, or $1$ for $i=1,\ldots, n{-}1$.
Consider a symmetry of $C'$ transposing the vectors $\bar v_k$
and $\bar v_{k+1}$ for $k=1,\ldots, n{-}2$. This transformation
is linear and its matrix coincides with the matrix of the
transposition of the vectors $\bar e_k$ and $\bar e_{k+1}$,
except the last column. The $n$-th column contains the
coordinates of the vector
$$
\begin{array}{c}
\frac{a_k-a_{k+1}}{2}(\bar e_k-\bar e_{k+1}) +\bar e_n.\\
\end{array}
$$
Since the transformation is lattice,
$$
\hbox{$a_k\equiv a_{k+1} (\MOD 2)$ for $k=1,\ldots, n{-}2$}.
$$
Since any $a_i$ is either zero or unit, the above imply
$a_1=a_2=\ldots =a_{n-1}$. If $a_1=0$, then the vector $\bar v_n$
is not of the unit lattice length, but the vector $\bar v_1$ is
of the unit length, so $C'$ is not lattice-regular. If $a_1=1$
then $C'$ coincides with $\{4,3^{n-2}\}^L_2$.

\vspace{2mm}

Suppose, that the faces of $(n{-}1)$-dimensional lattice-regular
cube $C$ are lattice-congruent to $\{4,3^{n-3}\}^L_2$. Then $C$
is lattice-congruent to the cube $C'$ generated by the origin and
the vectors $\bar v_i =\bar e_i$, for $i=1,\ldots, n-2$, $\bar
v_{n-1}=\bar e_1+\ldots+\bar e_{n-2}+2\bar e_{n-1}$, and  $\bar
v_n=a_1\bar e_1+\ldots+a_n\bar e_n$.

Consider a symmetry of $C'$ transposing the vectors $\bar
v_{n-1}$ and $\bar v_{n}$. This transformation is linear and its
matrix coincides with the matrix of identity transformation
except the last two columns. These columns contain the
coefficients of the following two vectors respectively:
$$
\begin{array}{c}
\sum\limits_{j=1}^{n}\frac{a_j}{2}\bar e_j-
\frac{1}{2}\sum\limits_{j=1}^{n-2}\bar e_j, \quad \hbox{and} \quad
\sum\limits_{j=1}^{n-2}\left(
\frac{a_{n-1}-2a_j-a_{n-1}a_j+2}{2a_n}\bar e_j\right)
+ \frac{4-a_{n-1}^2}{2a_n}\bar e_{n-1}-\frac{a_{n-1}}{2}\bar e_n.\\
\end{array}
$$
Since the described transformation is lattice, the integers
$a_{n-1}$ and $a_n$ are even, and $a_{n-2}$ is odd.

Consider now a symmetry of $C'$ transposing the vectors $\bar
v_{n-2}$ and $\bar v_{n-1}$. This transformation is linear and
its matrix coincides with the matrix of identity transformation
except the last three columns. These columns contain the
coefficients of the following three vectors respectively:
$$
\begin{array}{c}
\bar v_{n-1}, \quad -\sum\limits_{j=1}^{n-3}\bar e_j-\bar e_{n-1},
\quad \hbox{and} \quad
\frac{a_{n-1}-a_{n-2}}{a_n}\sum\limits_{j=1}^{n-3}\bar e_j+
2\frac{a_{n-1}-a_{n-2}}{a_n}\bar e_{n-1}+\bar e_n.\\
\end{array}
$$
Since the described transformation is lattice, the integer
$a_{n-1}-a_{n-2}$ is even, and thus $a_{n-2}$ is even. We come to
the contradiction with the divisibility of $a_{n-2}$ by 2. So $C$
is not lattice-regular.

\vspace{2mm}

Suppose, that the faces of $(n{-}1)$-dimensional lattice-regular
cube $C$ are lattice-congruent to $\{4,3^{n-3}\}^L_3$. Then $C$
is lattice-congruent to the cube $C'$ generated by the origin and
the vectors $\bar v_1 =\bar e_1$, $\bar v_i=\bar e_1 + 2\bar e_i$
for $i=2,\ldots,n{-}1$, and  $\bar v_n=a_1\bar e_1+\ldots+a_n\bar
e_n$. By Corollary~\ref{ind_dis_cube} we can choose $a_n$ equals
either $1$ or $2$. Then we choose $a_i$ being $0$, or $1$ for
$i=1,\ldots n{-}1$.

Consider a symmetry of $C'$ transposing the vectors $\bar
v_{n-1}$ and $\bar v_{n}$. This transformation is linear and its
matrix coincides with the matrix of identity transformation
except the last two columns. These columns contain the
coefficients of the following two vectors respectively:
$$
\begin{array}{c}
-\frac{1}{2}\bar e_1+\sum\limits_{j=1}^{n}\frac{a_j}{2}\bar e_j,
\quad \hbox{and} \quad \frac{(1-a_1)(2+a_{n-1})}{2a_n}\bar e_1-
\sum\limits_{j=2}^{n-2}\left( \frac{a_j(a_{n-1}+2)}{2a_n}\bar
e_j\right)+
\frac{4-a_{n-1}^2}{2a_n}\bar e_{n-1}-\frac{a_{n-1}}{2}\bar e_n.\\
\end{array}
$$
Since the described transformation is lattice, the integer $a_1$
is odd, and the integers $a_i$ for $i=1,\ldots n$ are even. Thus
$$
a_n=2, \qquad a_1=1, \quad \hbox{and} \quad a_2=\ldots=a_{n-1}=0.
$$
Therefore $C'$ coincides with $\{4,3^{n-2}\}^L_3$.

We have already studied all possible $n$-dimensional cases. This
proves the statement for the dimension $n$.

\vspace{2mm}

All statements of Proposition~\ref{regular_cubes} are proven. \qed

\section{Lattice-regular lattice generalized octahedra.}\label{octahedron}

In this section we describe all lattice-regular lattice
generalized octahedra for all integer dimensions greater than 2.

\begin{proposition}\label{regular_oct}
All elementary lattice-regular $n$-dimensional lattice
generalized octahedra for $n\ge 3$ are lattice-congruent to
$\{3^{n-2},4\}^L_1$, $\{3^{n-2},4\}^L_2$, or $\{3^{n-2},4\}^L_3$.
The generalized octahedra $\{3^{n-2},4\}^L_1$,
$\{3^{n-2},4\}^L_2$, and $\{3^{n-2},4\}^L_3$ are not
lattice-congruent to each other.
\end{proposition}

\begin{figure}[h]
$$\epsfbox{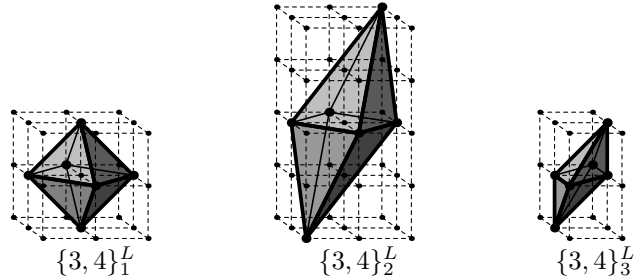}$$
\caption{Three-dimensional elementary lattice-regular convex
lattice tetrahedra.}\label{3doct}
\end{figure}

We show the (three-dimensional) octahedra $\{3,4\}^L_1$,
$\{3,4\}^L_2$, and $\{3,4\}^L_3$ on Figure~\ref{3doct}.

\begin{proof}

Consider an arbitrary elementary $n$-dimensional lattice-regular
generalized octahedron $P$. Let the vertices $V_1, \ldots,
V_{2n}$ of $P$ be enumerated in such a way that for any positive
integer $i\le n$ there exist a lattice symmetry exchanging $V_i$
and $V_{i+n}$ and preserving any other vertex. So,
$V_iV_jV_{i+n}V_{j+n}$ is a lattice-regular square for $i\ne j$.
Therefore, the midpoints of the segments $V_iV_{i+n}$ coincide
for $i=1,2,\ldots,n$. Denote the common midpoint of the segments
$V_iV_{i+n}$ by $A$.

\vspace{2mm}

Suppose the point $A$ is lattice. Consider the lattice cube with
the vertices $A\pm\bar{AV_1}\pm\bar{AV_2}\pm \ldots
\pm\bar{AV_n}$, we denote it by $C(P)$. The cube $C(P)$ is also
lattice-regular.

Note that the lattice-regular generalized octahedra $P'$ and $P''$
have lattice-congruent cubes $C(P')$ and $C(P'')$ iff $P'$ and
$P''$ are lattice-congruent.

Since $P$ is elementary, the segment $AV_1$ is of unit lattice
length. Therefore, the cube $C(P)$ is lattice-congruent to the
2-multiple of some $\{4,3^{n-2}\}^L_k$, for $k=1,2,3$. If $C(P)$
is lattice-congruent to the 2-multiple of $\{4,3^{n-2}\}^L_1$, or
to the 2-multiple of $\{4,3^{n-2}\}^L_2$, then $P$ is
lattice-congruent to $\{3^{n-2},4\}^L_1$, or to
$\{3^{n-2},4\}^L_2$, respectively. If $C(P)$ is lattice-congruent
to the 2-multiple of $\{4,3^{n-2}\}^L_3$, then $P$ is not
elementary.

\vspace{2mm}

Suppose now, the common midpoint $A$ of the diagonals is not
lattice. If the lattice length of $V_1V_{n+1}$ equals $2k{+}1$
for some positive $k$, then the generalized octahedron $P$ is not
elementary. Suppose the segment $V_1V_{n+1}$ is of unite lattice
length. Consider a 2-multiple to the polygon $P$ and denote in by
$2P$. Since the segment $V_1V_{n+1}$ is of unit lattice length,
the cube $C(2P)$ is the 2-multiple of  some $\{4,3^{n-2}\}^L_k$,
for $k=1,2,3$. If $C(2P)$ is lattice-congruent to the 2-multiple
of $\{4,3^{n-2}\}^L_1$, or of $\{4,3^{n-2}\}^L_2$, then $P$ is
not a lattice polytope. If $C(2P)$ is lattice-congruent to the
2-multiple of $\{4,3^{n-2}\}^L_3$, then $P$ is lattice-congruent
to $\{3^{n-2},4\}^L_3$.

\vspace{2mm}

The generalized octahedra $\{3^{n-2},4\}^L_1$,
$\{3^{n-2},4\}^L_2$, and $\{3^{n-2},4\}^L_3$ are lattice-regular,
since so are the cubes $C(\{3^{n-2},4\}^L_1)$,
$C(\{3^{n-2},4\}^L_2)$, and $C(2\{3^{n-2},4\}^L_3)$.

The generalized octahedra $\{3^{n-2},4\}^L_1$,
$\{3^{n-2},4\}^L_2$, and $\{3^{n-2},4\}^L_3$ are not
lattice-con\-gru\-ent to each-other, since the corresponding
elementary cubes $C(\{3^{n-2},4\}^L_1)$, $C(\{3^{n-2},4\}^L_2)$,
and $C(2\{3^{n-2},4\}^L_3)$ are not lattice-congruent.
\end{proof}

\section{Proof of Theorem~\ref{main_theorem}.}\label{3_4_dim}

In this section we obtain proof of Theorem~\ref{main_theorem} by
combining the results of propositions from the previous sections
and describing the remaining low-dimensional cases.

\vspace{2mm}

Consider any convex lattice-regular lattice polytope. Since it is
lattice-regular and convex it is combinatorically regular.
Therefore, by Theorem~B it is combinatorically isomorphic to one
of the Euclidean polytopes of Theorem~A. In Section~\ref{2d} we
gave the description of the two-dimensional case. In
Sections~\ref{symplex}, \ref{cube}, and~\ref{octahedron} we
studied the cases of lattice-regular polytopes combinatorically
isomorphic to regular symplices $(\{3^{n-1}\})$, regular cubes
$(\{4,3^{n-2}\})$, and regular generalized octahedra
$(\{3^{n-2},4\})$ respectively.

Now we will study the remaining special cases of three- and
four-dimensional regular polytopes.

\subsection{Three-dimensional icosahedra and dodecahedra.}
We have already classified all lattice-regular elementary
tetrahedra, cubes, and octahedra. There is no lattice-regular
dodecahedron, since there is no lattice-regular pentagon. There
is no lattice-regular icosahedron, since there is no
lattice-affine transformation with a fixed point of order 5.

So, the classification in the three-dimensional case is completed.

\subsection{Four-dimensional 24-sells, 120-sells, and 600-sells.}

{\it The case of 24-sell.} Suppose that $P$ is a lattice-regular
24-sell. It contains 16 vertices such that the subgroup of the
group of all lattice-symmetries of $P$ preserving these 16
vertices is isomorphic to the group of the symmetries of the
four-dimensional cube. So we can naturally define a
combinatorial-regular cube associated with this 16 vertices. Any
two-dimensional face of this cube is an Euclidean parallelogram,
since such face is the diagonal section containing four lattice
points of some lattice-regular octahedron. If all two-dimensional
faces are parallelograms, then all three-dimensional faces are
parallelepipeds and these 16 vertices are vertices of a
four-dimensional parallelepiped. Denote this parallelepiped by
$C$. Since all transformations of $C$ are lattice-affine, the
polytope $C$ is a lattice-regular cube. (Note that for any
$24$-sell there exist exactly three such (distinct) cubes).

Suppose that $C$ is a lattice-regular cube generated by the origin
and some vectors $\bar v_i$ for $i=1,2,3,4$. Consider also the
coordinates $(*,*,*,*)_v$ corresponding to this basis. Let the
point $(a_1,a_2,a_3,a_4)_v$ of $P$ connected by edges with the
vertices of $C$ is in the plane with the unit last coordinate.
Then, the point $(2{-}a_1,a_2,a_3,2{-}a_4)_v$ is also a vertex.
Thus, the point $(a_1,a_2,a_3,a_4{-}1)$ is also a vertex. Note
that the points $(a_1,a_2,a_3,a_4)_v$ and
$(a_1,a_2,a_3,a_4{-}2)_v$ are symmetric about the center of the
cube: $(1/2,1/2,1/2,1/2)_v$. So $a_4{-}1/2=1/2-a_4{-}2$. Thus,
$a_4=3/2$. Similar calculations show that $a_1=a_2=a_3=1/2$.
Therefore, the eight points do not contained in $C$ coincide with
the following points:
$$
O+1/2(\bar v_1{+}\bar v_2{+}\bar v_3{+}\bar v_4)\pm\bar v_i,\quad
\hbox{for $i=1,2,3,4$.}
$$

Consider a lattice-regular primitive 24-sell $P$ and $C$ one of
the corresponding cubes. Since the edges of $C$ are the edges of
$P$, the cube $C$ is also primitive. Let us study all three
possible cases of lattice-affine types of $C$.

Suppose $C$ coincides with $\{4,3,3\}^L_1$. Then the remaining
points
$$
O+1/2(\bar e_1{+}\bar e_2{+}\bar e_3{+}\bar e_4)\pm\bar e_i,\quad
\hbox{for $i=1,2,3,4$}
$$
are not lattice. Therefore, the case of $\{4,3,3\}^L_1$ is
impossible.

In the case of $\{4,3,3\}^L_2$ and $\{4,3,3\}^L_3$ all the
vertices are lattice and in our notation coincide with
$\{3,4,3\}^L_1$ and $\{3,4,3\}^L_2$ respectively. Straightforward
calculations shows, that both resulting 24-sells are
lattice-regular.

{\it The case of 120-sell.} The 120-sell is non-realizable as an
lattice-regular lattice polytope, since its two-dimensional faces
should be lattice-regular lattice pentagons. By the above,
lattice-regular lattice pentagons are not realizable.

{\it The case of 600-sell.} Consider an arbitrary polytope with
topological structure of 600-sell having one vertex at the origin
$O$. Let $OV$ be some edge of this 600-sell. The group of
symmetries of an abstract 600-sell with fixed vertex $O$ is
isomorphic to the group of symmetries of an abstract icosahedron.
The group of symmetries of an abstract 600-sell with fixed
vertices $O$ and $V_1$ is isomorphic to the group of symmetries
of an abstract pentagon. So, there exist a symmetry $A$ of the
600-sell of order 5 preserving the vertices $O$ and $V_1$. If the
polytope $P$ is lattice-regular, then this symmetry is
lattice-linear. Since $A^5$ is the identity transformation and
the space is four-dimensional, the characteristic polynomial of
$A$ in the variable $x$ is either $x{-}1$ or
$x^4{+}x^3{+}x^2{+}x{+}1$. Since $A(\bar{OV_1})=\bar{OV_1}$, the
characteristic polynomial of $A$  is divisible by $x{-}1$. If the
characteristic polynomial is $x{-}1$, then the operator $A$ is
the identity operator of order 1 and not of order 5. Therefore,
there is no lattice-regular lattice polytope with the
combinatorial structure of the 600-sell.

\subsection{Conclusion of the proof of Theorem~\ref{main_theorem}.}

We have studied all possible combinatorical cases of
lattice-regular polytopes. The proof of
Theorem~\ref{main_theorem} is completed. \qed


\begin{thebibliography}{99}
\bibitem{Arn2}
V.~I.~Arnold, {\it Continued fractions}, M.: Moscow Center of
Continuous Mathematical Education, 2002.
\bibitem{Arn5}
V.~I.~Arnold, {\it Statistics of integer convex polygons}, Func.
an. and appl., v.~14(1980), n.~2, pp.~1--3.
\bibitem{Bar1}
I.~B\'ar\'any J.-M.~Kantor, {\it On the number of lattice free
polytopes}, European J. Combin. v.~21(1), 2000, pp.~103--110.
\bibitem{Bar2}
I.~B\'ar\'any A.~M.~Vershik, {\it On the number of convex lattice
polytopes}, Geom. Funct. Anal. v.~2(4), 1992, pp.~381--393.
\bibitem{Cox}
H.S.M.~Coxeter, {\it The complete enumeration of finite groups of
the form $R^2_i=(R_iR_j)^{k_{ij}}=1$}, J. London Math. Soc, v.~10
(ser.~1), No 37, 1935, pp.~21--25.
\bibitem{Dan}
L.~Danzer, E.~Schulte, {\it Regul\"are inzidenzkomplexe I}, Geom.
Dedicata, v.~13(1982), pp.~295--308.
\bibitem{Ful}
W.~Fulton, {\it Introduction to Toric Varieties}, Annals of
Mathematics Studies; Princeton University Press, v.~131(1993),
\bibitem{JMK1}
J.-M.~Kantor, K.~S.~Sarakaria {\it On primitive subdivisions of
an elementary tetrahedron}, Pacific J. of Math. v.~211(1), 2003,
pp.~123--155.
\bibitem{KarPyr}
O.~N.~Karpenkov, {\it Two-dimensional faces of multidimensional
continued fractions and pyramids on integer lattices}, preprint n
2005-26, Cahiers du CEREMADE, Unite Mixte de Recherche du
C.N.R.S. N7534 (2005),
http://www.ceremade.dauphine.fr/preprint/CMD/2005-26.pdf.
\bibitem{Kh1}
A.~G.~Khovanskii, A.~Pukhlikov, {\it Finitely additive measures
of virtual polytopes}, Algebra and Analysis, v.~4(2), 1992,
pp.~161--185.
\bibitem{Kor2}
E.~I.~Korkina, {\it Two-dimensional continued fractions. The
simplest examples}, Proceedings of V.~A.~Steklov Math. Ins., v.
209(1995), pp. 243--166.
\bibitem{LacBook}
G.~Lachaud, {\it Voiles et polyhedres de Klein}, Act. Sci. Ind.,
176 pp., Hermann, 2002.
\bibitem{McM1}
P.~McMullen, {\it Combinatorially regular polytopes},
Mathematika, v.~14(1967), pp~142--150.
\bibitem{McM2}
P.~McMullen, {\it Affinely and projectively regular polytopes},
J. London Math. Soc, v.~43(1968), pp~755--757.
\bibitem{McM3}
P.~McMullen, E.~Schulte, {\it Abstract Regular Polytopes},
Series: Encyclopedia of Mathematics and its Applications,
No.~92(2002), 566 p.
\bibitem{Mou2}
J.-O.~Moussafir, {\it Voiles et Poly\'edres de Klein: Geometrie,
Algorithmes et Statistiques},
docteur en sciences th\'ese, Universit\'e Paris IX - Dauphine, (2000),\\
see also at http://www.ceremade.dauphine.fr/\~{}msfr/.
\bibitem{Oda}
T.~Oda, {\it Convex bodies and Algebraic Geometry, An
Introduction to the Theory of Toric Varieties}, Springer-Verlag,
Survey in Mathemayics, 15(1988).
\bibitem{Rez}
B.~Reznick, {\it Lattice point symplices}, Discrete Math.,
v.~60(1986), pp.~219--242.
\bibitem{Sch}
L.~Schl\"afli, {\it Theorie der vielfachen Kontinuit\"at},
written 1850-1852; Z\"urcher und Furrer, Z\"urich 1901;
Denks\-chriften der Schweizerischen naturforschenden Gesellschaft
v.~38(1901), pp.~1--237; reprinted in {\it Ludwig Schl\"afli
1814-1895, Gesammelte Mathematische Abhandlungen}, Vol.~I,
Birkh\"auser, Basel 1950, pp.~167--387.
\bibitem{Whi}
G.~K.~White, {\it Lattice tetrahedra}, Canadian J. of Math.
16(1964), pp.~389--396.
\end{thebibliography}
\end{document}